\documentclass[12pt]{amsart}
\usepackage{amssymb}
\usepackage{amsmath}
\usepackage{mathtools}
\usepackage{a4wide}
\usepackage{color}
 \usepackage{boxedminipage}
\usepackage{enumerate, graphicx}
\usepackage{todonotes}
\usepackage{subfig}
\usepackage{wrapfig}
\usepackage{cite}
\usepackage{caption}
\usepackage{float}
\usepackage[colorlinks,citecolor=blue,linkcolor=blue]{hyperref}
\usepackage[capitalise,nameinlink]{cleveref}
\usepackage{url}

\def\Var{{\mathfrak V}}
\def\Tscr{{\mathcal T}}

  \makeatletter
  \CheckCommand*\refstepcounter[1]{\stepcounter{#1}%
      \protected@edef\@currentlabel
       {\csname p@#1\endcsname\csname the#1\endcsname}%
  }
  \renewcommand*\refstepcounter[1]{\stepcounter{#1}%
    \protected@edef\@currentlabel
      {\csname p@#1\expandafter\endcsname\csname the#1\endcsname}%
  }
  \def\labelformat#1{\expandafter\def\csname p@#1\endcsname##1}
  \DeclareRobustCommand\Ref[1]{\protected@edef\@tempa{\ref{#1}}%
     \expandafter\MakeUppercase\@tempa
  }
  \makeatother

  \makeatletter
  \newcommand{\numberlike}[2]{%
     \expandafter\def\csname c@#1\endcsname{%
         \expandafter\csname c@#2\endcsname}%
  }
  \makeatother

  \def\DefaultNumberTheoremWithin{section}

  \theoremstyle{plain}
  
     \numberwithin{Lemma}{\DefaultNumberTheoremWithin}
     \labelformat{Lemma}{Lemma~#1}
  
     \numberwithin{Claim}{\DefaultNumberTheoremWithin}
     \numberlike{Claim}{Lemma}
     \labelformat{Claim}{Claim~#1}

  \newtheorem{Theorem}{Theorem}
     \numberwithin{Theorem}{\DefaultNumberTheoremWithin}
     \numberlike{Theorem}{Lemma}
     \labelformat{Theorem}{Theorem~#1}
  
     \numberwithin{Corollary}{\DefaultNumberTheoremWithin}
     \numberlike{Corollary}{Lemma}
     \labelformat{Corollary}{Corollary~#1}
  
     \numberwithin{Proposition}{\DefaultNumberTheoremWithin}
     \numberlike{Proposition}{Lemma}
     \labelformat{Proposition}{Proposition~#1}
  
     \numberwithin{Conjecture}{\DefaultNumberTheoremWithin}
     \numberlike{Conjecture}{Lemma}
     \labelformat{Conjecture}{Conjecture~#1}
  
     \numberwithin{Situation}{\DefaultNumberTheoremWithin}
     \numberlike{Situation}{Lemma}
     \labelformat{Situation}{Situation~#1}
 
     \numberwithin{Note}{\DefaultNumberTheoremWithin}
     \numberlike{Note}{Lemma}
     \labelformat{Note}{Note~#1}
     
  \theoremstyle{definition}
  \newtheorem{Definition}{Definition}
     \numberwithin{Definition}{\DefaultNumberTheoremWithin}
     \numberlike{Definition}{Lemma}
     \labelformat{Definition}{Definition~#1}

  \theoremstyle{definition}
  
     \numberwithin{Question}{\DefaultNumberTheoremWithin}
     \numberlike{Question}{Lemma}
     \labelformat{Question}{Question~#1}

  \theoremstyle{definition}
  
     \numberwithin{Problem}{\DefaultNumberTheoremWithin}
     \numberlike{Problem}{Lemma}
     \labelformat{Problem}{Problem~#1}

     \theoremstyle{remark} 
     \numberwithin{Remark}{\DefaultNumberTheoremWithin}
     \numberlike{Remark}{Lemma}
     \labelformat{Remark}{Remark~#1}
  \theoremstyle{remark}

     \numberwithin{Example}{\DefaultNumberTheoremWithin}
     \numberlike{Example}{Lemma}
     \labelformat{Example}{Example~#1}
  
     \labelformat{Case}{Case~#1}
     \numberwithin{Case}{Lemma}
  
     \labelformat{Step}{Step~#1}
     \numberwithin{Step}{Lemma}

\labelformat{section}{Section~#1}
  \labelformat{subsection}{Section~#1}
  \labelformat{subsubsection}{Section~#1}

\title{Exploring the Varchenko Determinant of Partial Cubes}
\author{Winfried Hochst\"attler}
    \address{FernUniversit\"at in Hagen \\ 
          Fakult\"at f\"ur Mathematik und Informatik \\
          58084 Hagen\\
          Germany}
     \email{winfried.hochstaettler@fernuni-hagen.de}
     
\author{Sophia Keip}
    \address{FernUniversit\"at in Hagen \\ 
          Fakult\"at f\"ur Mathematik und Informatik \\
          58084 Hagen\\
          Germany}
     \email{sophia.keip@fernuni-hagen.de}

\author{Birol Yazici}
    \address{FernUniversit\"at in Hagen \\ 
          Fakult\"at f\"ur Mathematik und Informatik \\
          58084 Hagen\\
          Germany}
     \email{yazici.birol@gmail.com}

\begin{document}
\begin{abstract}
The Varchenko matrix is known to have a well-structured determinant for complexes of oriented matroids (COMs). COMs can be characterized as partial cubes that do not have certain forbidden pc-minors. In this work, we generalize the Varchenko matrix and its determinant to partial cubes. We identify examples of partial cubes whose Varchenko determinants lack a clean factorization, as well as those that exhibit such a structure. These findings open the door for further research into the properties and potential characterizations of partial cubes with well-behaved Varchenko determinants.
\end{abstract}

\maketitle
\section{Introduction}
In 1993 Varchenko introduced a symmetric bilinear form on the regions of affine hyperplane arrangements \cite{varchenko1993bilinear}. The rows and columns of the matrix representing that bilinear form are indexed by the regions of the arrangement $\mathcal{H}$
and an entry belonging to two regions $Q_i$ and $Q_j$ is calculated as $\prod_{e \in S(Q_i, Q_j)} w_e$, where the $w_e$ represent weights assigned to the hyperplanes $H_e \in \mathcal{H}$. The set $S(Q_i, Q_j)$ consists of the hyperplanes that must be crossed along a shortest path between $Q_i$ and $Q_j$ and is called separator. Varchenko showed, that the determinant of that matrix has a very elegant factorization, where each factor looks like
\begin{align}\label{eq:factor}
    (1-(\prod_{i \in I}w_i)^2)^b.
\end{align}
$I$ is the index set belonging to a non empty intersection of hyperplanes in $\mathcal{H}$ and $b$ is a non negative integer. Since then this result has been generalized to cones \cite{gente2013varchenko}, oriented matroids (OMs) \cite{hochstattler2019varchenko, randriamaro2020varchenko}, topoplane arrangements \cite{randriamaro2022varchenko} and complexes of oriented matroids (COMs) \cite{hochstattler2025signed}.
COMs can be characterized as partial cubes which do not contain certain minors. The Varchenko matrix extends naturally to partial cubes. In \cite{hochstattler2025signed}, the authors posed the question whether there are classes of partial cubes beyond COMs that also have a well-structured Varchenko determinant. To address this question, we computed the Varchenko determinant for several characterizing excluded minors. While the smallest minor yielded an unstructured determinant, suggesting that only COMs might exhibit a clean factorization, further exploration disproved this: we also found a minor with a nicely factorized determinant. We will begin with an introduction to COMs, followed by their connection to partial cubes. Next, we revisit the determinant formula for COMs and present examples of non-COM partial cubes with both structured and unstructured determinants. This work serves as a foundation for further exploration of the Varchenko matrix in the context of partial cubes.

\section{Complexes of Oriented Matroids}
COMs were introduced as a generalization of lopsided sets and oriented matroids \cite{bandelt2018coms}. A COM is defined by a finite ground set \( E \) and a collection of sign vectors \( \mathcal{L} \subseteq \{0, +, -\}^{|E|} \). Before introducing the axioms that characterize a COM, we first define the separator and the composition of two sign vectors.

\begin{Definition} Let $E$ be a finite set and $X,Y \in \{0,+,-\}^{|E|}$. We call
\begin{align*}
    S(X,Y) = \{e \in E: X_e = -Y_e \neq 0\}. 
\end{align*}
the \emph{separator} of $X$ and $Y$. The \emph{composition} of $X$ and $Y$ is defined by
\begin{align*}
(X \circ Y)_e = \begin{cases}
                        X_e & \text{ if } X_e \neq 0,\\
                        Y_e & \text{ if } X_e = 0\\
                       \end{cases} 
                    \forall e\in E.
\end{align*} 
\end{Definition}

With these definitions in place, we can now proceed to define COMs.

\begin{Definition}[Complex of Oriented Matroids (COM)] Let $E$ be a finite set and $\mathcal{L} \subseteq \{0,+,-\}^{|E|}$. The pair $\mathcal{M}=(E,\mathcal{L})$ is called a COM, if $\mathcal{L}$ satisfies 
\begin{itemize}
    \item[(FS)]  Face Symmetry 
    \begin{align*}
\forall X,Y \in \mathcal{L}: X \circ (-Y) \in \mathcal{L}.
\end{align*}
\end{itemize}
and
\begin{itemize}
\item[(SE)] Strong Elimination 
\begin{align*}
&\forall X,Y \in \mathcal{L}\, \forall e \in S(X,Y)\, \exists Z \in \mathcal{L}: \\
&Z_e=0 \text{ and }\forall f \in E \setminus S(X,Y): Z_f = (X \circ Y)_f.
\end{align*} 
\end{itemize}
The elements of $\mathcal{L}$ are called \emph{covectors}.
\end{Definition}

\begin{figure}
    \centering
    \includegraphics[scale=0.9]{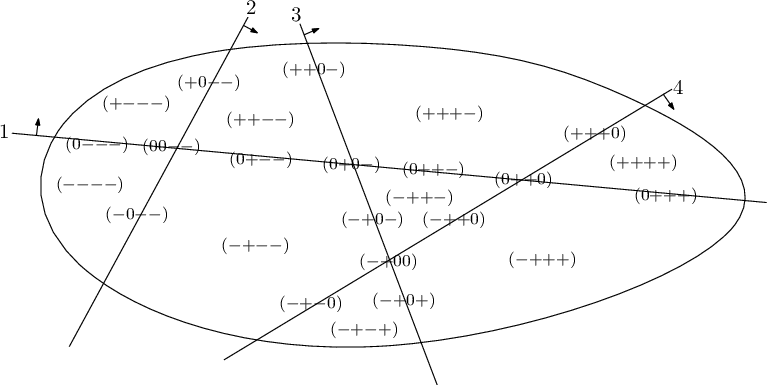}
    \caption{Example of a realizable COM induced by the intersection of four hyperplanes with an open convex set. The positive side of each hyperplane is indicated by an arrow. Each sign vector encodes whether the corresponding cell lies on the positive side, the negative side, or directly on the hyperplane.}
    \label{fig:COM}
\end{figure}

A special case of COMs are the so-called realizable COMs, which provide a concrete and intuitive example of this structure. These arise from the intersection of an oriented affine hyperplane arrangement with an open convex set. To construct a realizable COM, we examine all cells of the hyperplane arrangement within the open convex set. For each cell, we record its position relative to the hyperplanes in a sign vector, indicating whether it lies on the positive side, the negative side, or directly on the hyperplane. The resulting collection of sign vectors forms the covector set of a COM. In this particular case, the separators correspond to the hyperplanes that distinguish one cell from another. An example of a realizable COM is illustrated in \cref{fig:COM}. 

Let $X \in \{0,+,-\}^{|E|}$. We define the \emph{support} of $X$, denoted $\underline{X}$, as $\underline{X} = \{e \in E: X_e \neq  0\}$. The zero set $E\backslash \underline{X}$ of a covector $X$ is denoted by $z(X)$. 
Covectors of a COM $\mathcal{M} = (E, \mathcal{L})$ with maximal support are called \emph{topes}. We denote the set of all topes of a COM by $\mathcal{T}$. In the following, we assume that $\mathcal{M}$ is \emph{simple}. This means that for all \(e \in E\), we have $\{X_e \mid X \in \mathcal{L}\} = \{+, -, 0\}$, and for all distinct $e, f \in E$, it holds that $\{X_e X_f \mid X \in \mathcal{L}\} = \{+, -, 0\}$. 
Under this assumption, the topes are precisely the covectors that do not contain any zero entries. Let us now define the tope graph of a COM:

\begin{Definition}[Tope Graph of a COM]
Let $\mathcal{M} = (E, \mathcal{L})$ be a COM and $\mathcal{T}$ its tope set. $G = (V,E)$ is called the tope graph of $\mathcal{M}$ if $V = \mathcal{T}$ and $(X,Y) \in E$ if $|S(X,Y)|=1$ for all $X,Y \in V$.
\end{Definition}

\begin{figure}
    \centering
    \includegraphics[scale=0.9]{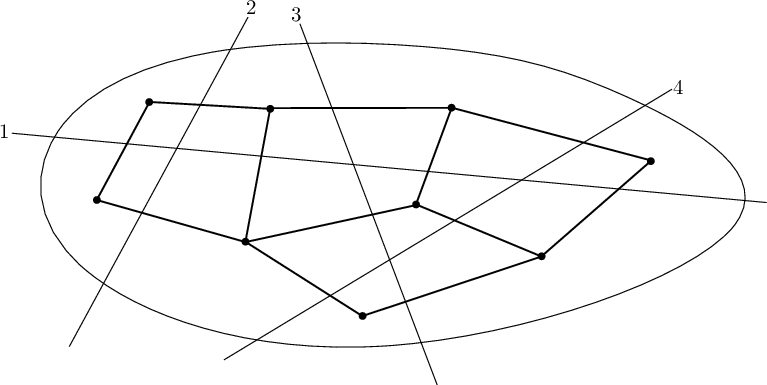}
    \caption{Tope graph of the realizable COM from \cref{fig:COM}.}
    \label{fig:topegraph}
\end{figure}

Using this definition, the tope graph of the realizable COM shown in \cref{fig:COM} can be constructed by assigning a vertex to each full-dimensional cell of the hyperplane arrangement and connecting two vertices if their corresponding cells are separated by exactly one hyperplane. This is illustrated in \cref{fig:topegraph}.

In \cite[Corollary 4.10]{knauer2020tope} the authors showed the following.

\begin{Theorem} Let $\mathcal{M} = (E,\mathcal{L})$ be a simple COM, $\mathcal{T}$ its tope set and $G = (V,E)$ its tope graph. $\mathcal{M}$ is uniquely detemined by $\mathcal{T}$ and up to reorientation by $G$. 
\end{Theorem}

In the next chapter we will get to know partial cubes, which are closely connected to the tope graphs of COMs.

\section{Partial Cubes}
A hypercube graph is a graph with $2^n$ vertices, corresponding to all binary strings of length $n$. Two vertices are connected by an edge if the corresponding binary strings only differ by one element. Therefore every vertex has $n$ adjacent neighbors. We denote the $n$-dimensional hypercube by $Q_n$. A graph $G = (V,E)$ is called partial cube if it is a isometric subgraph of the hypercube. This means that the distance~$d$ - that is, the length of a shortest path between two vertices - is preserved: $d_G(u,v) = d_{Q_n}(u,v)$ for all $u,v \in V$.
 Knauer and Marc showed the following in \cite{bandelt2018coms}.
\begin{Theorem}
  Let $\mathcal{M}=(E,\mathcal{L})$ be a simple COM. Then the tope graph of $\mathcal{M}$ is a partial cube.  
\end{Theorem}
Moreover, they showed that COMs can be defined as partial cubes that do not have certain minors. To understand this kind of minors, we need to define equivalence classes on the edges of a partial cube. For that we need the Djokovi\'c-Winkler-Relation: 
\begin{Definition}[Djokovi\'c-Winkler-Relation \cite{djokovic1973distance,winkler1984isometric}]
Let $G=(V,E)$ be a partial cube and $(u,v) \in E$. Let
\begin{align*}
W(u,v) =  \{x \in V : d(x,u) < d(x,v)\}
\end{align*}
and $e = (u,v),f = (x,y) \in E$. The \emph{Djokovi\'c-Winkler-Relation} $\theta$ is defined by
\begin{align*}
    e\, \theta f  \Leftrightarrow  x\in W(u,v) \land y \in W(v,u).
\end{align*}
\end{Definition}

If $G$ is a partial cube, the Djokovi\'c-Winkler-Relation defines equivalence classes on the edges of the graph. We call these the \emph{color classes} $\mathcal{E}$ of the partial cube $G$. The sets $\mathcal{E}^+ = W(u,v)$ and $\mathcal{E}^-=W(v,u)$ are called \emph{complementary halfspaces} of $G$. 

\begin{figure}
       \centering
       \includegraphics[scale=0.9]{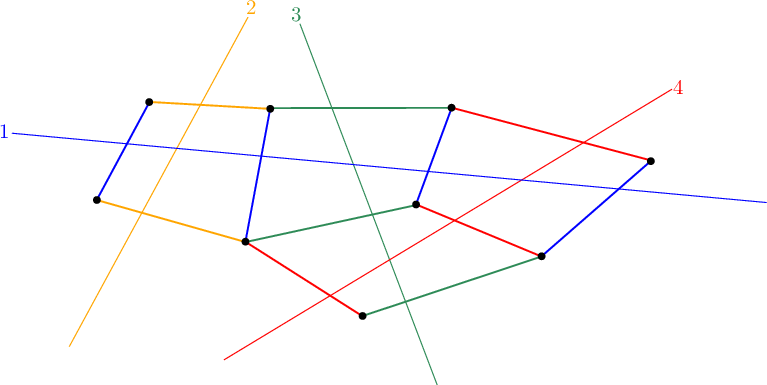}\\
       \caption{The color classes of the realizable COM in \cref{fig:colorclasses}}
       \label{fig:colorclasses}
   \end{figure}

The color classes of a partial cube representing the tope graph of a realizable COM can be determined straightforwardly: all edges corresponding to the same hyperplane belong to the same color class. This is visualized in \cref{fig:colorclasses} for our COM from \cref{fig:COM}. In this case, the complementary halfspaces correspond to the sets of vertices on either side of the hyperplane. 

More generally, the color classes of a partial cube group together edges corresponding to changes in the same bit of the binary strings in the underlying hypercube. One of the complementary halfspaces consist of all vertices that have a value of $1$ at a specific bit and the other of the vertices that have a value of $0$.

\begin{figure}
       \centering
      \hfill A) \includegraphics[scale=0.9]{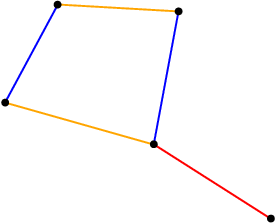}\hfill
      B) \includegraphics[scale=0.9]{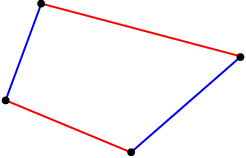}\hfill
      C) \includegraphics[scale=0.9]{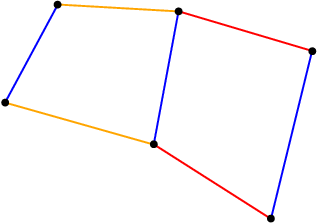} \hfill
       \caption{Pc-minors of the realizable COM from \cref{fig:COM}. In A) the graph is restricted to the left side of the green color class, in B) to the right side, and in C) the entire green color class is contracted.}
       \label{fig:pcminors}
   \end{figure}

A \emph{partial cube minor} (pc-minor) is obtained by either contracting an entire color class or restricting the graph to one of its half spaces. It represents a specific type of classical graph minor. For the realizable COM shown in \cref{fig:COM}, the different ways for deriving a pc-minor are illustrated in \cref{fig:pcminors}.

Let us now define the class of forbidden pc-minors, that characterize COMs. 
\begin{Definition}[Forbidden pc-minors]
Let $Q_n$, $n\geq 4$ be the hypercube and fix a vertex $v$. Delete the antipode $-v$ of $v$ from $Q_n$ together with 
\begin{itemize}
    \item[-] either exactly one neighbor of $v$, denoted by $Q_n^{-*}$,
    \item[-] or $v$ and at least one of its neighbors, denoted by $Q_n^{--}(m),\,n\geq m>0$.
\end{itemize}

We set $\mathcal{Q}^- = \{Q_n^{-*}, Q_n^{--}(m)\,|\,4 \leq n, \,0 < m \leq n\}$ and call this class of partial cubes the \emph{forbidden pc-minors}.    
\end{Definition}

\begin{figure}
       \centering
    \begin{minipage}[t]{0.23\textwidth}
        \centering
        $Q_4$ \\
        \includegraphics[scale=0.6]{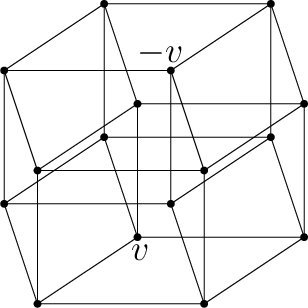}
    \end{minipage}
    \hfill
    \begin{minipage}[t]{0.23\textwidth}
        \centering
        $Q_4^{-*}$ \\
        \includegraphics[scale=0.6]{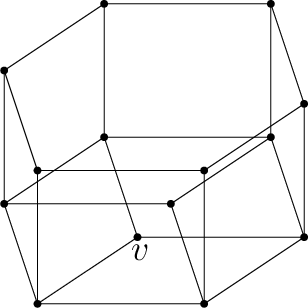}
    \end{minipage}
    \hfill
    \begin{minipage}[t]{0.23\textwidth}
        \centering
        $Q_4^{--}(1)$ \\
        \includegraphics[scale=0.6]{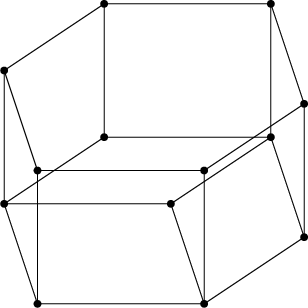}
    \end{minipage}
    \hfill
    \begin{minipage}[t]{0.23\textwidth}
        \centering
        $Q_4^{--}(4)$ \\
        \includegraphics[scale=0.6]{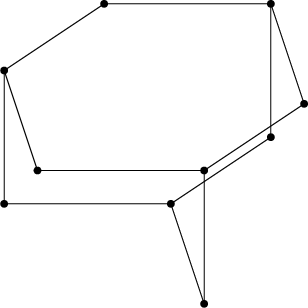}
    \end{minipage}
       \caption{
        Illustration of forbidden pc-minors: On the left, the vertex \( v \) and its antipode \( -v \) are highlighted within the hypercube \( Q_4 \). To obtain each forbidden pc-minor, \( -v \) must be removed. Additionally, for \( Q_4^{--} \), exactly one neighbor of \( v \) is deleted. For \( Q_4^{--}(1) \) and \( Q_4^{--}(4) \), \( v \) is removed along with 1 or 4 of its neighbors, respectively.
        }

       \label{fig:forbiddenpcminors}
   \end{figure}

For $n=4$ some of the forbidden pc-minors are demonstrated in \cref{fig:forbiddenpcminors}.
Knauer and Marc \cite[Theorem 1.1]{knauer2020tope} proved the following:
\begin{Theorem}
Let $G$ be a partial cube. $G$ is the tope graph of a COM if and only if $G$ has no pc-minor in the set $\mathcal{Q}^-$.
\end{Theorem}
Since the tope graph of a COM determines it completely up to orientation, the theorem means that COMs can be characterized as partial cubes that do not contain a forbidden pc-minor.

To summarize, forbidden pc-minors provide a clear criterion for distinguishing COMs among partial cubes. As we aim to explore the Varchenko determinant of partial cubes beyond COMs, these minors serve as natural candidates for investigation, offering minimal examples of non-COM structures.

\section{The Varchenko Determinant}

We start by defining the Varchenko matrix of a COM:

\begin{Definition}[Varchenko Matrix of a COM]\label{def:signed} Let $\mathcal{M}=(E,\mathcal{L})$ be a COM and $\mathcal{T}$ its tope set. The Varchenko matrix $\Var$ of a COM is a $\# \Tscr \times \# \Tscr$-Matrix over 
\begin{align*}
\mathbb{K}[x_{e}\,|\, e\in E].
\end{align*}
Its rows and columns are indexed by the topes $\Tscr$ in a fixed linear order. For $P,Q \in \Tscr$
\begin{align*}
\mathfrak{V}_{P,Q} = \prod_{e \in S(P,Q)} x_{e}.
\end{align*}
\end{Definition}

The Varchenko matrix of the realizable COM $\mathcal{M}$ from \cref{fig:COM} is

\begin{align*}
\Var(\mathcal{M}) &= \begin{bsmallmatrix}
1 & x_2 & x_2 x_3 & x_2 x_3 x_4 & x_1 & x_1 x_2 & x_1 x_2 x_3 & x_1 x_2 x_3 x_4 & x_1 x_2 x_4 \\
x_2 & 1 & x_3 & x_3 x_4 & x_1 x_2 & x_1  & x_1 x_3 & x_1 x_3 x_4 & x_1 x_4\\
x_2 x_3 & x_3 & 1 & x_4 & x_1 x_2 x_3 & x_1 x_3 & x_1 & x_1 x_4 & x_1 x_3 x_4\\
x_2 x_3 x_4 & x_3 x_4 & x_4 & 1 & x_1 x_2 x_3 x_4 & x_1 x_3 x_4 & x_1 x_4 & x_1 & x_1 x_3 \\
x_1 & x_1 x_2 & x_1 x_2 x_3 & x_1 x_2 x_3 x_4 & 1 & x_2 & x_2 x_3 & x_2 x_3 x_4 & x_2 x_4\\
x_1 x_2 & x_1 & x_1 x_3 & x_1 x_3 x_4 & x_2 & 1 & x_3 & x_3 x_4 & x_4\\ 
x_1 x_2 x_3 & x_1 x_3 & x_1 & x_1 x_4 & x_2 x_3 & x_3 & 1 & x_4 & x_3 x_4\\
x_1 x_2 x_3 x_4 & x_1 x_3 x_4 & x_1 x_4 & x_1 & x_2 x_3 x_4 & x_3 x_4 & x_4 & 1 & x_3\\
x_1 x_2 x_4 & x_1 x_4 & x_1 x_3 x_4 & x_1 x_3 & x_2 x_4 & x_4 & x_3 x_4 & x_3 & 1
\end{bsmallmatrix}
\end{align*}

In \cite{hochstattler2025signed}, the authors showed that the determinant of the Varchenko matrix of a COM has a structured factorization formula:
\begin{Theorem}\label{thm:varchenko}
   Let $\Var$ be the Varchenko matrix of the COM $\mathcal{M} = (E,\mathcal{L})$.
   Then
   \begin{align*}
     \det (\Var) = \prod_{Y \in \mathcal{L}} (1-(\prod_{e \in z(Y)} x_{e})^2)^{b_Y}. 
   \end{align*}
   where $b_Y$ are nonnegative integers.
\end{Theorem}

We observe that, instead of indexing the variables in the monomials of the individual factors by intersections of hyperplane arrangements, we now have zero sets of covectors. This correspondence becomes especially clear in the realizable case. For the Varchenko matrix of the COM in \cref{fig:COM} we computed the determinant and get 
\begin{align*}
    \text{Det}(\Var(\mathcal{M})) = (1-x_1^2)^4(1-x_2^2)^2(1-x_3^2)^3(1-x_4^2)^3
\end{align*}
which is in line with \ref{thm:varchenko}.

The Varchenko matrix can be easily generalized to partial cubes by substituting the role of the topes with vertices and defining the separator $S(u,v)$ of two vertices $u$ and $v$ by the color classes we have to cross on a shortest path between two vertices. 

\begin{Definition}[Varchenko Matrix of a Partial Cube]\label{def:varpartial} Let $G = (V,E)$ be a partial cube and let $\mathcal{E}$ be its color classes. The Varchenko matrix $\Var$ of a partial cube is a $|V| \times |V|$-Matrix over 
\begin{align*}
\mathbb{K}[x_{e}\,|\, e \in \mathcal{E}].
\end{align*}
Its rows and columns are indexed by $V$ in a fixed linear order. For $u,v \in V$ the entries of the Varchenko matrix are defined as
\begin{align*}
\mathfrak{V}_{u,v} = \prod_{e \in S(u,v)} x_{e}.
\end{align*}
\end{Definition}

\begin{figure}
\includegraphics[scale = 0.6]{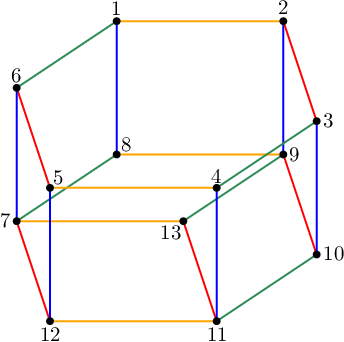} \hspace{2cm}
\includegraphics[scale = 0.6]{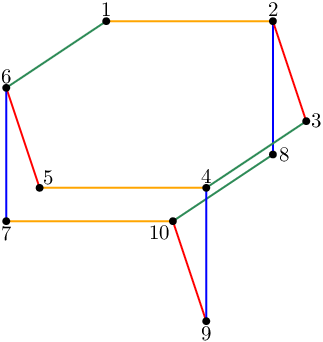} 
\caption{The color classes of $Q_4^{--}(1)$ and $Q_4^{--}(4)$. The vertex numbering corresponds to the linear order that we assume when indexing the rows and columns of the Varchenko matrix.}\label{fig:Q4colorclasses}
\end{figure}

Let us demonstrate this on the forbidden pc-minor $Q_4^{--}(4)$. The color classes of this partial cube are illustrated in \cref{fig:Q4colorclasses} together with the assumed linear order of the vertices. We will use $x_1$ for the yellow color class, $x_2$ for the red one and $x_3$, $x_4$ for green and blue respectively.
We get the following Varchenko matrix for $Q_4^{--}(4)$:
\begin{align*}
\Var(Q_4^{--}(4)) &= \begin{bsmallmatrix}
1 & x_1 & x_1 x_2 & x_1 x_2 x_3 & x_2 x_3 & x_3 & x_3 x_4 & x_1 x_4 & x_1 x_2 x_3 x_4 & x_1 x_3 x_4 \\
x_1 & 1 & x_2 & x_2 x_3 & x_1 x_2 x_3 & x_1 x_3 & x_1 x_3 x_4 & x_4 & x_2 x_3 x_4 & x_3 x_4 \\
x_1 x_2 & x_2 & 1 & x_3 & x_1 x_3 & x_1 x_2 x_3 & x_1 x_2 x_3 x_4 & x_2 x_4 & x_3 x_4 & x_2 x_3 x_4 \\
x_1 x_2 x_3 & x_2 x_3 & x_3 & 1 & x_1& x_1 x_2 & x_1 x_2 x_4 & x_2 x_3 x_4 & x_4 & x_2 x_4 \\
x_2 x_3 & x_1 x_2 x_3 & x_1 x_3 & x_1 & 1 & x_2 & x_2 x_4 & x_1 x_2 x_3 x_4 & x_1 x_4 & x_1 x_2 x_4 \\
x_3 & x_1 x_3 & x_1 x_2 x_3 & x_1 x_2 & x_2 & 1 & x_4 & x_1 x_3 x_4 & x_1 x_2 x_4 & x_1 x_4 \\
x_3 x_4 & x_1 x_3 x_4 & x_1 x_2 x_3 x_4 & x_1 x_2 x_4 & x_2 x_4 & x_4 & 1 & x_1 x_3 & x_1 x_2 & x_1 \\
x_1 x_4 & x_4 & x_2 x_4 & x_2 x_3 x_4 & x_1 x_2 x_3 x_4 & x_1 x_3 x_4 & x_1 x_3 & 1 & x_2 x_3 & x_3 \\
x_1 x_2 x_3 x_4 & x_2 x_3 x_4 & x_3 x_4 & x_4 & x_1 x_4 & x_1 x_2 x_4 & x_1 x_2 & x_2 x_3 & 1 & x_2 \\
x_1 x_3 x_4 & x_3 x_4 & x_2 x_3 x_4 & x_2 x_4 & x_1 x_2 x_4 & x_1 x_4 & x_1 & x_3 & x_2 & 1 \\
\end{bsmallmatrix}
\end{align*}
The question of whether the generalized Varchenko determinant retains a clean factorization, i.e. the factors look like in (\ref{eq:factor}), beyond COMs has been answered negatively in \cite{hochstattler2025signed} by examining $Q_4^{--}(4)$, as
\begin{align*}
\text{Det}(\Var(Q_4^{--}(4))) = & (1-x_1^2)^3 (1-x_2^2)^3 (1-x_3^2)^3 (1-x_4^2)^3 \dots\\
&\dots(3x_1^2x_2^2x_3^2x_4^2 - x_1^2x_2^2x_3^2 - x_1^2x_2^2x_4^2 - x_1^2x_3^2x_4^2 - x_2^2x_3^2x_4^2 + 1)
\end{align*}
where the last factor is obviously not in the form of (\ref{eq:factor}). This leads to another question: Does the Varchenko matrix of a partial cube factorize nicely if and only if it is the tope graph of a COM? This question can also be answered negatively by looking at the Varchenko determinant of $Q_4^{--}(1)$. The color classes of that partial cube together with the linear order of the vertices are also illustrated in \cref{fig:colorclasses}. 
\begin{align*}
&\Var(Q_4^{--}(1)) = \\
&\hspace{-1cm}
\begin{bsmallmatrix}
1 & x_1 & x_1 x_2 & x_1 x_2 x_3 & x_2 x_3 & x_3 & x_3 x_4 & x_4 & x_1 x_4 & x_1 x_2 x_4 & x_1 x_2 x_3 x_4& x_2 x_3 x_4 & x_1 x_3 x_4 \\
x_1 & 1 & x_2 & x_2 x_3 & x_1 x_2 x_3 & x_1 x_3 & x_1 x_3 x_4 & x_1 x_4 & x_4 & x_2 x_4 & x_2 x_3 x_4& x_1 x_2 x_3 x_4 & x_3 x_4 \\
x_1 x_2 & x_2 & 1 & x_3 & x_1 x_3 & x_1 x_2 x_3 & x_1 x_2 x_3 x_4 & x_1 x_2 x_4 & x_2 x_4 & x_4 & x_3 x_4& x_1 x_3 x_4 & x_2 x_3 x_4 \\
x_1 x_2 x_3 & x_2 x_3 & x_3 & 1 & x_1 & x_1 x_2 & x_1 x_2 x_4 & x_1 x_2 x_3 x_4 & x_2 x_3 x_4 & x_3 x_4 & x_4& x_1 x_4 & x_2 x_4 \\
x_2 x_3 & x_1 x_2 x_3 & x_1 x_3 & x_1 & 1 & x_2 & x_2 x_4 & x_2 x_3 x_4 & x_1 x_2 x_3 x_4 & x_1 x_3 x_4 & x_1 x_4& x_4 & x_1 x_2 x_4 \\
x_3 & x_1 x_3 & x_1 x_2 x_3 & x_1 x_2 & x_2 & 1 & x_4 & x_3 x_4 & x_1 x_3 x_4 & x_1 x_2 x_3 x_4 & x_1 x_2 x_4& x_2 x_4 & x_1 x_4 \\
x_3 x_4 & x_1 x_3 x_4 & x_1 x_2 x_3 x_4 & x_1 x_2 x_4 & x_2 x_4 & x_4 & 1 & x_3 & x_1 x_3 & x_1 x_2 x_3 & x_1 x_2 & x_2 & x_1 \\
x_4 & x_1 x_4 & x_1 x_2 x_4 & x_1 x_2 x_3 x_4 & x_2 x_3 x_4 & x_3 x_4 & x_3 & 1 & x_1 & x_1 x_2 & x_1 x_2 x_3 & x_2 x_3 & x_1 x_3 \\
x_1 x_4 & x_4 & x_2 x_4 & x_2 x_3 x_4 & x_1 x_2 x_3 x_4 & x_1 x_3 x_4 & x_1 x_3 & x_1 & 1 & x_2 & x_2 x_3& x_1 x_2 x_3 & x_3 \\
x_1 x_2 x_4 & x_2 x_4 & x_4 & x_3 x_4 & x_1 x_3 x_4 & x_1 x_2 x_3 x_4 & x_1 x_2 x_3 & x_1 x_2 & x_2 & 1 & x_3 & x_1 x_3 & x_2 x_3 \\
x_1 x_2 x_3 x_4 & x_2 x_3 x_4 & x_3 x_4 & x_4 & x_1 x_4 & x_1 x_2 x_4 & x_1 x_2 & x_1 x_2 x_3 & x_2 x_3 & x_3 & 1& x_1 & x_2 \\
x_2 x_3 x_4 & x_1 x_2 x_3 x_4 & x_1 x_3 x_4 & x_1 x_4 & x_4 & x_2 x_4 & x_2 & x_2 x_3 & x_1 x_2 x_3 & x_1 x_3 & x_1 & 1 & x_1 x_2  \\
x_1 x_3 x_4 & x_3 x_4 & x_2 x_3 x_4 & x_2 x_4 & x_1 x_2 x_4 & x_1 x_4 & x_1 & x_1 x_3 & x_3 & x_2 x_3 & x_2 & x_1 x_2 & 1 \\
\end{bsmallmatrix}
\end{align*}

For this matrix we get
\begin{align*}
\text{Det}(\Var(Q_4^{--}(1))) = (1-x_1^2)^5(1-x_2^2)^5(1-x_3^2)^5(1-x_4^2)^6(1-(x_1x_2x_3)^2)
\end{align*}
where each factor is in the form of (\ref{eq:factor}). This observation raises the concluding question of whether partial cubes with a well-structured Varchenko determinant can be characterized by other distinct properties.

\section{Conclusion}
In this paper, we generalized the Varchenko matrix to partial cubes and observed its determinant. Although the Varchenko matrix admits a clean factorization when the partial cube is the tope graph of a COM, this property does not hold in general, as illustrated by a known example of a minimal non-COM partial cube. We answered the question whether the Varchenko determinant has a nice factorization \emph{if and only if} the partial cube is the tope graph of a COM negatively: by examining the determinants of other partial cubes that are forbidden pc-minors for COMs, we found one that still admits a clean factorization. This raises the question of whether partial cubes with a nice Varchenko determinant factorization can be characterized in another way. If this question can be answered, it would be very interesting to know to what structural property the indices of the monomials in the factorization belong.

\section*{Appendix: Computations}
To feed the intuition for further research on Varchenko determinants of partial cubes, we provide a summary of all the determinants we tested below. We used Python and Maple for our computations. In order to check if a graph is a partial cube and determine its color classes, we used the work of David Eppstein \cite{eppstein2007recognizing} together with the corresponding code \url{https://ics.uci.edu/~eppstein/PADS/}. The Python code that generates the Varchenko matrix of a partial cube can be found here: 
\begin{center}
    \url{https://github.com/SophiaKeip/VarchenkoMatrixPartialCubes}.
\end{center}
 Note that the color classes and vertices in the code are numbered differently from those in the examples above. The numbering in the examples was chosen for illustrative purposes.
\vspace{3mm}

Varchenko determinants for $n = 4$:
\allowdisplaybreaks
\begin{align*}
\Var(Q^{*-}_4) =& 
(1-x_1^2)^6(1-x_2^2)^6(1-x_3^2)^6(1-x_4^2)^6 (1-(x_1x_3x_4)^2)\\[10pt]
\Var(Q^{--}_4(1)) =& 
(1-x_1^2)^6(1-x_2^2)^5(1-x_3^2)^5(1-x_4^2)^5(1-(x_2x_3x_4)^2)\\[10pt]
\Var(Q^{--}_4(2)) =&
(1-x_1^2)^5(1-x_2^2)^5(1-x_3^2)^4(1-x_4^2)^4 
\\&\left( x_1^2x_2^2x_3^2x_4^2 - x_1^2x_3^2x_4^2 - x_2^2x_3^2x_4^2 + 1 \right)\\[10pt]
\Var(Q^{--}_4(3)) =& 
(1-x_1^2)^4(1-x_2^2)^4(1-x_3^2)^4(1-x_4^2)^3 
\\&\left( 2x_1^2x_2^2x_3^2x_4^2 - x_1^2x_2^2x_4^2 - x_1^2x_3^2x_4^2 - x_2^2x_3^2x_4^2 + 1 \right)\\[10pt]
\Var(Q^{--}_4(4)) =& 
(1-x_1^2)^3(1-x_2^2)^3(1-x_3^2)^3(1-x_4^2)^3 
\\&\left( 3x_1^2x_2^2x_3^2x_4^2 - x_1^2x_2^2x_3^2 - x_1^2x_2^2x_4^2 - x_1^2x_3^2x_4^2 - x_2^2x_3^2x_4^2 + 1 \right)
\end{align*}

Varchenko determinants for $n = 5$:
\allowdisplaybreaks
\begin{align*}
\Var(Q^{*-}_5) = 
&(1-x_1^2)^{14}(1-x_2^2)^{14}(1-x_3^2)^{14}(1-x_4^2)^{14}(1-x_5^2)^{14}(1-(x_2x_3x_4x_5)^2)\\[10pt]
\Var(Q^{--}_5(1)) = 
&(1-x_1^2)^{13}(1-x_2^2)^{13}(1-x_3^2)^{14}(1-x_4^2)^{13}(1-x_5^2)^{13}(1-(x_1x_2x_4x_5)^2)\\[10pt]
\Var(Q^{--}_5(2)) = 
&(1-x_1^2)^{13}(1-x_2^2)^{12}(1-x_3^2)^{13}(1-x_4^2)^{12}(1-x_5^2)^{12}\\&
\left( x_1^2x_2^2x_3^2x_4^2x_5^2 - x_1^2x_2^2x_4^2x_5^2 - x_2^2x_3^2x_4^2x_5^2 + 1 \right)\\[10pt]
\Var(Q^{--}_5(3)) = 
&(1-x_1^2)^{12}(1-x_2^2)^{12}(1-x_3^2)^{12}(1-x_4^2)^{11}(1-x_5^2)^{11}\\&
\left( 2x_1^2x_2^2x_3^2x_4^2x_5^2 - x_1^2x_2^2x_4^2x_5^2 - x_1^2x_3^2x_4^2x_5^2 - x_2^2x_3^2x_4^2x_5^2 + 1 \right)\\[10pt]
\Var(Q^{--}_5(4)) = 
&(1-x_1^2)^{11}(1-x_2^2)^{11}(1-x_3^2)^{11}(1-x_4^2)^{11}(1-x_5^2)^{10}\\&
\left( 3x_1^2x_2^2x_3^2x_4^2x_5^2 - x_1^2x_2^2x_3^2x_5^2 - x_1^2x_2^2x_4^2x_5^2 - x_1^2x_3^2x_4^2x_5^2 - x_2^2x_3^2x_4^2x_5^2 + 1 \right)\\[10pt]
\Var(Q^{--}_5(5)) = 
&(1-x_1^2)^{10}(1-x_2^2)^{10}(1-x_3^2)^{10}(1-x_4^2)^{10}(1-x_5^2)^{10}\\&
( 4x_1^2x_2^2x_3^2x_4^2x_5^2 - x_1^2x_2^2x_3^2x_4^2 - x_1^2x_2^2x_3^2x_5^2 - x_1^2x_2^2x_4^2x_5^2- x_1^2x_3^2x_4^2x_5^2 - x_2^2x_3^2x_4^2x_5^2 + 1 ) 
\end{align*}

\bibliographystyle{plain} 
\bibliography{pcbib}

\end{document}